\def\hybrid{\topmargin 0pt      \oddsidemargin 0pt
        \headheight 0pt \headsep 0pt
        \textwidth 160true mm       
        \textheight 231true mm         
        \marginparwidth 0.0in
        \parskip 0pt plus 1pt   \jot = 1.5ex}
\catcode`\@=11
\def\marginnote#1{}

\newcount\hour
\newcount\minute
\newtoks\amorpm
\hour=\time\divide\hour by60
\minute=\time{\multiply\hour by60 \global\advance\minute by-\hour}
\edef\standardtime{{\ifnum\hour<12 \global\amorpm={am}%
        \else\global\amorpm={pm}\advance\hour by-12 \fi
        \ifnum\hour=0 \hour=12 \fi
        \number\hour:\ifnum\minute<10 0\fi\number\minute\the\amorpm}}
\edef\militarytime{\number\hour:\ifnum\minute<10 0\fi\number\minute}

\def\draftlabel#1{{\@bsphack\if@filesw {\let\thepage\relax
   \xdef\@gtempa{\write\@auxout{\string
      \newlabel{#1}{{\@currentlabel}{\thepage}}}}}\@gtempa
   \if@nobreak \ifvmode\nobreak\fi\fi\fi\@esphack}
        \gdef\@eqnlabel{#1}}
\def\@eqnlabel{}
\def\@vacuum{}
\def\draftmarginnote#1{\marginpar{\raggedright\scriptsize\tt#1}}

\def\draft{\oddsidemargin -.5truein
        \def\@oddfoot{\sl preliminary draft \hfil
        \rm\thepage\hfil\sl\today\quad\militarytime}
        \let\@evenfoot\@oddfoot \overfullrule 3pt
        \let\label=\draftlabel
        \let\marginnote=\draftmarginnote
   \def\@eqnnum{(\theequation)\rlap{\kern\marginparsep\tt\@eqnlabel}%
\global\let\@eqnlabel\@vacuum}  }

\relax
\hybrid
\draft

\documentclass[10pt]{article}
\usepackage{amsmath}
\usepackage{amsfonts}
\def\ba{\begin{eqnarray}}
\def\ea{\end{eqnarray}}
\def\y{\hat{y}}
\def\lb{\label}
\def\Tr{{\rm Tr}}
\def\be{\begin{equation}}
\def\ee{\end{equation}}
\def\qed{\rule{5pt}{5pt}}

\newcommand{\bea}{\begin{eqnarray}}
\newcommand{\eea}{\end{eqnarray}}

\newtheorem{thm}{Theorem}

\newtheorem{prop}{Proposition}[section]

\def\Z{\mathbb{Z}}

\def\S{\mathbb{S}}

\def\C{\mathbb{C}}

\begin{document}

{}

\vspace{2cm}

\begin{center}
{\Large \bf Bethe subalgebras in Hecke algebra and Gaudin models}
\footnote{The work of A.P.Isaev was supported by the grant RFBR 11-01-00980-a
and grant Higher School of Economics No.11-09-0038.}
\end{center}

\vspace{1cm}

\begin{center}
\large{A.P.~Isaev$^*$ and A.N.~Kirillov$^{**}$}
\end{center}

\begin{center}
$^*$ Bogoliubov Laboratory of Theoretical Physics, JINR, \\
 141980, Dubna, Moscow region, \\
 and ITPM, M.V.Lomonosov Moscow State University, Russia \\
 E-mail: isaevap@theor.jinr.ru
\end{center}

\begin{center}
$^{**}$ Research Institute of Mathematical Sciences, RIMS, \\
Kyoto University, Sakyo-ku, 606-8502, Japan \\
\end{center}


{\bf Abstract.} The generating function for elements of the Bethe subalgebra
of Hecke algebra is constructed as Sklyanin's transfer-matrix operator for Hecke chain.
We show that in a special classical limit $q \to 1$ the
Hamiltonians of the Gaudin model can be derived from the transfer-matrix operator of Hecke chain. We consruct a non-local  analogue of the Gaudin
Hamiltonians for the case of Hecke algebras.

\newpage

\section{Introduction}

The Gaudin models were firstly introduced by M.Gaudin in \cite{Gaud}.
  These models were also investigated as limiting cases of integrable
quantum inhomogeneous
$su(2)$-chains in \cite{Skly}.
 Here we use
 an algebraic approach and obtain Gaudin's Hamiltonians from the
transfer-matrix operator for open inhomogeneous chain models
which formulated in terms of generators of affine Hecke algebra
$\hat{H}_{M+1}(q)$. In our chain model an inhomogeneity
appears  (as well as in \cite{Skly}) as different shifts in spectral parameters related to different sites
of the Hecke chain. The Gaudin  Hamiltonians are obtained from the generating
function which defines a Bethe subalgebra in the Hecke algebra
$\hat{H}_{M+1}(q),$  by taking a special ``classical limit'' $q \to 1$.

The Bethe subalgebras in the group ring of symmetric groups have been studied,
for example, in \cite{MTV}, \cite{K}. In the present paper we construct a
{\bf lift} of the Bethe subalgebras studied in the papers mentioned above, to
 the cases of the Hecke and affine Hecke algebras. Our
construction of the Bethe subalgebras is based on some special properties of
the {\bf trace} maps \cite{Jo},\cite{IsOg1},\cite{Isa07} in the {\it tower} of the (affine)
Hecke algebras, see  Section~3, formulae $(15).$  ~Non formally speaking, the
main idea behind our construction, is to define  a set of ``baxterized''
Jucys--Murphy elements in the (affine) Hecke algebra. To realize this idea we
treat the Jucys--Murphy elements in the (affine) Hecke algebra as a ``classical
limit'' of the canonical free abelian subgroup in the (affine) braid group,
see Section~2.

The plan of the paper is as follows. In Section~2 we review some basic facts
about braid and affine braid groups we need, namely, definitions and the
construction of the maximal free abelian subgroups in these groups. \\
Sections~3,4
contain our main results, namely, the construction of Bethe's subalgebras in
the Hecke and affine Hecke algebras. In particular, Theorem 1
in Section~4 describes the Hecke version of the Gaudin  Hamiltonians. We treat the Bethe subalgebras obtained
 as a ``baxterization'' of the canonical free abelian subgroup in the
corresponding braid group. In other words, we introduce a spectral parameter
dependences in definition of the  Jucys--Murphy elements keeping
the commutativity property of the deformed elements. A similar construction
can be done for the Birman--Murakami-Wenzl algebras, cyclotomic Hecke
algebras and some other quotients of braid groups.

In Section~5 we study the classical and Yangian limits of the Bethe subalgebras
 in the Hecke and affine Hecke algebras correspondingly.
 
 \vspace{0.2cm}
 
 We thank S.Krivonos for valuable discussions.   
 
\section {Braid group}

Denote by ~$\mathbb{S}_{n}$~\underline {\bf the symmetric group}
 on $n$ letters, and by~ $s_{i}$ ~the simple transposition
$(i,i+1)$ for $1 \le i \le n-1.$

The well-known Moore--Coxeter presentation of the symmetric group has the form
$$\langle s_{1},\ldots,s_{n-1} ~~\vert ~~s_{i}^2=1,~~s_{i}s_{i+1}s_{i}=s_{i+1}s_{i}s_{i+1}
,~s_{i}s_{j}=s_{j}s_{i},~~~if~~~|i-j| \ge 2 \rangle.$$

Transpositions $s_{ij}:=s_i~s_{i+1} \cdots s_{j-2}~s_{j-1}~
s_{j-2} \cdots s_{i+1}~s_{i},$~$1 \le i < j < j \le n,$  satisfy the following
set of (defining) relations:

$s_{ij}^{2}=1,$~~$s_{ij}~s_{kl}=s_{kl}~s_{ij},$~ if~ $\{i,j\} \bigcap \{k,l \}=
\emptyset,$~~

$s_{ij}~s_{ik}=s_{jk}~s_{ij}=s_{ik}~s_{jk},$~~$s_{ik}~s_{ij}=s_{ij}~s_{jk}=
s_{jk}~s_{ik},$~~$i < j < k.$

{\underline {{\bf The Artin braid group on $n$ strands $B_{n}$}} is defined by
generators $\sigma_{1},\ldots,\sigma_{n-1}$ and relations
\begin{equation}
\sigma_{i}\sigma_{i+1}\sigma_{i}=\sigma_{i+1}\sigma_{i}\sigma_{i+1},~~
1 \le i \le n-2, ~~~\sigma_{i}\sigma_{j}=\sigma_{j}\sigma_{i} ~~if ~~|i-j|
\ge 2.
\end{equation}

\begin{prop}
 Let us introduce elements
$$D_{i,j}:=\sigma_{j-1}\sigma_{j-2}\cdots\sigma_{i+1}~
\sigma_{i}^2~\sigma_{i+1}\cdots\sigma_{j-2}\sigma_{j-1},$$
$$F_{i,j}:=
\sigma_{n-j} \sigma_{n-j+1} \ldots \sigma_{n-i-1} ~\sigma_{n-i}^{2}
~\sigma_{n-i-1} \ldots \sigma_{n-j+1} \sigma_{n-j},$$
where $1 \le i < j \le n.$

For example, \\
$D_{i,i+1}=\sigma_{i}^2,~~D_{i,i+2}=
\sigma_{i+1}\sigma_{i}^2\sigma_{i+1},~~~ F_{i,i+1}=\sigma_{n-i}^{2},$~~
$F_{i,i+2}=\sigma_{n-i-1}~\sigma_{n-i}^2~\sigma_{n-i-1},$~
 and so on. \\
\underline{Then}

$\bullet$~~~For each~$j=3, \ldots,n,$~the element $D_{1,j}$~commutes with
~$\sigma_{1},\ldots,\sigma_{j-2}.$

$\bullet$ The elements $D_{i,i+1},D_{i,i+2},\ldots,D_{i,n}$ (resp.
$F_{i,i+1}, F_{i,i+2},\ldots,F_{i,n}$) ~$1 \le i \le n-1,$~~generate a free
abelian subgroup in $B_{n}.$

$\bullet$ The elements $D_{1,2},D_{2,3},\ldots,D_{1,n}$ (resp.
$F_{1,2}, F_{2,3},\ldots,F_{1,n}$) ~~generate a {\it maximal}  free
abelian subgroup in $B_{n}.$

$\bullet$ ~~~If $n \ge 3,$ the element
$$ \prod_{2 \le j \le n}D_{1,j} = \prod_{2 \le j \le n}~F_{1,j}=
(\sigma_{1}\cdots\sigma_{n-1})^n $$
generates the center of the braid group $B_{n}.$

$\bullet$~~~  $D_{i,j}D_{i,j+1}D_{j,j+1}=D_{j,j+1}D_{i,j+1}D_{i,j},~$
if ~$i < j.$
\end{prop}
$\bullet$~~Consider the elements $s:=\sigma_1~\sigma_2~\sigma_1,$~
$t:= \sigma_1~\sigma_2$ in the braid group $B_3.$ Then $s^2=t^3$ and the
element $c:=s^2$ generates the center of the group $B_3.$ Moreover,
$$ B_3 / \langle c \rangle \cong PSL_2(\Z),~~
B_3 / \langle c^2 \rangle \cong SL_2(\Z).$$

{\underline {\bf The affine Artin braid group $B_n^{aff}$}},
is an
extension of the Artin Braid group on $n$ strands $B_n$ by the element $\tau$
subject to the set of crossing relations
$$ \sigma_1~\tau~\sigma_1~\tau=\tau~\sigma_1~\tau~\sigma_1,
~~\sigma_i~\tau=\tau~\sigma_i~~for~~2 \le i \le n-1.$$
\begin{prop} The elements
$${\hat D}_{1}:= \tau,~ {\hat D}_{j}= \sigma_{j-1}~ {\hat D}_{j-1}~ \sigma_{j-1},~~~2 \le j \le n,$$
generate a free abelian subgroup in $B_n^{aff}.$
\end{prop}
Therefore, for a unital commutative algebra $F,$ any quotient
$F[B_n] / J$ of the group algebra $F[B_n]$ of the braid group $B_n$ (resp.
a quotient $F[B_n^{aff}] / I$ of the affine braid group $B_n^{aff}$ group
algebra $F[B_n^{aff}]$) by a two-sided ideal $J \subset F[B_n]$
(resp $I \subset F[B_n^{aff}]$) contains distinguish  commutative subalgebra
generated by the images of elements $D_{1,2},\ldots,D_{1,n-1}$~
(resp. ${\hat D}_1, \ldots,{\hat D}_n$). It
 is well-known that the Hecke and  affine Hecke algebras are certain
quotients of the braid and affine braid groups correspondingly, see the next
Section for details. In these cases the images of elements $D_{1,2},\ldots,
D_{1,n-1}$ and those ${\hat D}_1,\ldots, {\hat D}_n$ coincide with the Jucys--
Murphy elements in the Hecke and affine Hecke algebras correspondingly. Our
main objective of the next Section is to construct an analogue of the Bethe
subalgebras in the affine Hecke algebras.

 \section{Bethe subalgebras for affine Hecke algebra}

The {\underline {\bf Hecke algebra $H_{M+1}(q)$}}
(see, e.g.,  \cite{ChPr} and \cite{Jo}) is generated by invertible
elements $T_i$ $(i=1, \dots, M)$ subject to the set of relations:
\be
\label{braidg}
T_i \, T_{i+1} \, T_i =
T_{i+1} \, T_i \,  T_{i+1} \; , \;\;\;
T_{i} \,  T_{j} = T_{j} \,  T_{i} \;\;
{\rm for} \;\; i \neq j \pm 1 \; .
\ee
\be
\label{ahecke}
T^2_i - 1 = (q-q^{-1}) \,  T_i \; ,
\ee
Let $x$ be a spectral parameter. We define
\underline{Baxterized elements}
\be
\lb{baxtH}
T_i(x) : =  T_i - x T_i^{-1} = (1-x)T_i + \lambda x \in H_{M+1}(q) \; ,
\ee
which in view of (\ref{braidg}) and (\ref{ahecke}) solve the Yang-Baxter
equation
\be
\lb{ybeH}
T_i(x) \, T_{i-1}(xz) \, T_i(z) =
T_{i-1}(z) \, T_i(xz) \, T_{i-1}(x)  \; ,
\ee
and satisfy relations
\be
\lb{bax1}
T_i(x)T_i(z)=\lambda T_i(x z) + (1-x)(1-z) \; ,
\ee
\be
\lb{bax2}
\begin{array}{c}
T_i(1) = \lambda \; , \;\;\;
T_i(x) = \frac{(1-x)}{(1-z)}T_i(z) +\lambda \frac{(x-z)}{(1-z)} \; ,
0\end{array}
\ee
where $\lambda := (q -q^{-1})$. Note that from (\ref{bax1}) for Baxterized
elements we have the condition
$$
T_i(x) \, T_i(x^{-1})= (q^{-1} x - q) (q^{-1} x^{-1} - q) \; ,
$$
which can be written as unitarity condition $\widetilde{T}_i(x) \widetilde{T}_i(x^{-1}) =1$
for modified baxterized elements
 \be
\lb{bax3}
\widetilde{T}_i(x) = \frac{1}{(q- q^{-1} x)} \, T_i(x) \; .
\ee

The {\underline {\bf affine Hecke algebra $\hat{H}_{M+1}(q)$}}~~ (see, e.g.,
Chapter 12.3 in \cite{ChPr} and \cite{Kirill}) is an extension of the Hecke
algebra$H_{M+1}(q)$ by additional affine elements $\y_k$ $(k=1,\dots,M+1)$
subject to relations:
\be
\lb{afheck}
\y_{k+1} = T_k \, \y_k \, T_k \; , \;\;\; \y_k \, \y_j = \y_j \, \y_k \; , \;\;\;
 \y_j \, T_i  =  T_i \, \y_j \;\; (j \neq i,i+1) \; .
\ee
The elements $\{ \y_k \}$ form a commutative subalgebra in $\hat{H}_{M+1}$,
while the symmetric functions in $\y_k$ form the center in $\hat{H}_{M+1}$.
The Jucys--Murphy elements $\{ \y_k \}$ coincide with the images of elements
${\hat D}_1,\ldots, {\hat D}_n$ considered in previous Section.
Here and below we omit the dependence on $q$ in the notations $H_{M+1}(q)$ and
$\hat{H}_{M+1}(q)$ of the Hecke  algebras.

The {\underline {\bf Ariki-Koike algebra \cite{ArKo},\cite{BrMa}  ${\cal H}_{M+1}(q,Q_1,...,Q_m)$}}~~
is the quotient of the affine Hecke algebra $\hat{H}_{M+1}$ by the characteristic identity
 \be
\lb{charH}
 (\y_1 - Q_1) \cdots (\y_1 - Q_m) = 0 \; ,
\ee
where $Q_1, \dots, Q_m$ are parameters.

\vspace{0.2cm}

{\bf Definition~3.1}~~ Let  $\vec{\xi}_{(n)}=(\xi_1,\dots,\xi_{n})$ be $n$
parameters and $y_1(x) \in \hat{H}_{M+1},$ define the elements
\be
\lb{xxz55}
\begin{array}{c}
y_{n}(x;\vec{\xi}_{(n-1)}) =
T_{n-1}(\frac{x}{\xi_{n-1}}) \cdots T_2(\frac{x}{\xi_2}) T_1(\frac{x}{\xi_1}) \, y_1(x) \,
T_{1}(x \xi_1) T_2(x \xi_2) \cdots T_{n-1}(x \xi_{n-1}) = \\ [0.2cm]
= T_{n-1}(\frac{x}{\xi_{n-1}}) y_{n-1}(x;\vec{\xi}_{(n-2)}) T_{n-1}(x \xi_{n-1}) \; ,
\end{array}
\ee
which we call as ``baxterized'' Jucys--Murphy elements.

{\bf Proposition~3.1}~~ {\it Assume that the element $y_1(x) \in \hat{H}_{M+1}$ in
 (\ref{xxz55}) is any {\it local}
(i.e., $[y_1(x), T_k]=0,$ $\forall k > 1$) solution of the reflection equation
\be
\lb{reflH1}
T_1 \left(x /z \right) \, y_{1}(x) \, T_1(x \, z) \, y_{1}(z)) =
y_{1}(z) \, T_1(x \, z) \, y_{1}(x) \,  T_1 \left( x /z \right) \; ,
 \ee

Then the elements  (\ref{xxz55}) satisfy  the
\underline{ reflection equation}
 \be
 \lb{reflH}
T_n \left(x /z \right) \, y_{n}(x;\vec{\xi}_{(n-1)}) \, T_n(x \, z) \, y_{n}(z;\vec{\xi}_{(n-1)}) =
y_{n}(z;\vec{\xi}_{(n-1)}) \, T_n(x \, z) \, y_{n}(x;\vec{\xi}_{(n-1)}) \,  T_n \left( x /z \right) \; ,
 \ee
}
{\bf Proof}~~The case $n=1$ of the equation (\ref{reflH})  corresponds to our
assumption that $y_1(x)$ satisfies the equation (\ref{reflH1}). The general
case follows by induction using the definition (\ref{xxz55}) of elements
$y_{n}(x;\vec{\xi}_{(n-1)}).$ \\
\qed

For example, in the case of the affine Hecke algebra, one can use the local solution (see \cite{IsOg1}):
\be
\lb{soluH}
 y_1(x)= { {\hat{y}}_1 -\xi x \over {\hat{y}_1 - \xi x^{-1}}} ,
 \ee
 where $\xi$ is a parameter.
 In the case of the Ariki-Koike algebra this rational solution is represented
 in the polynomial form by writing the characteristic identity (\ref{charH}) as
 $$
 { 1 \over {\hat{y}_1 - \xi x^{-1}}} = v_1 \, \hat{y}_1^{m-1} + v_2 \, \hat{y}_1^{m-2} + \dots +
 v_{m-1} \, \hat{y}_1 + v_{m} \; ,
 $$
 where $v_1,\dots,v_m$ are functions of $\xi, x, Q_1, \dots, Q_m$.
\vspace{0.3cm}

Consider the following inclusions of the subalgebras
$\hat{H}_{1} \subset \hat{H}_{2} \subset \dots \subset \hat{H}_{M+1}$:
$$
\{\y_1; T_1, \dots ,T_{n-1}\} \in \hat{H}_{n} \subset \hat{H}_{n+1} \ni
\{\y_1; T_1, \dots ,T_{n-1},T_n \} \; .
$$
Define for the algebra $\hat{H}_{M+1}$ linear mappings
$$
{\rm Tr}_{(n+1)}: \;\; \hat{H}_{n+1} \to \hat{H}_{n} \; , \;\;\;\; (n=1,2,\dots,M) \; ,
$$
such that for all $X,X' \in
\hat{H}_{n}$ and $Y \in \hat{H}_{n+1}$ we have
\be
\label{map}
\begin{array}{c}
 {\rm Tr}_{(n+1)} ( T_n^{\pm 1} \cdot X \cdot T_n^{\mp 1}) = {\rm Tr}_{(n)} (X)   \, , \;\;
{\rm Tr}_{(n+1)} ( X \cdot Y \cdot X' ) = X \cdot {\rm Tr}_{(n+1)}(Y) \cdot X' \;\;  \, , \\[0.1cm]
{\rm Tr}_{(n)}  {\rm Tr}_{(n+1)} ( T_n \cdot Y ) = {\rm Tr}_{(n)}  {\rm Tr}_{(n+1)}
( Y \cdot T_n)  \; , \;\;\; \\[0.1cm]
{\rm Tr}_{(n+1)} (T_n) =  1 \; ,
{\rm Tr}_{(1)} (y_1^k)= D^{(k)} \; , \;\; {\rm Tr}_{(n+1)} ( X ) = D^{(0)} \, X
 \; ,
\end{array}
\ee
where $k \in \mathbb{Z}$ and $D^{(k)} \in \mathbb{C}\backslash \{0\}$  are constants. Note that
$D^{(0)}$ is independent of $n$ and all
$D^{(k)}$ can be considered as central elements for certain central extension Ext$(\hat{H}_{M+1})$
of $\hat{H}_{M+1}$. The elements
$D^{(k)}$ generate an abelian subalgebra (we denote this subalgebra $\hat{H}_{0}$) in Ext$(\hat{H}_{M+1})$.

Using the properties (\ref{map}) of the map ${\rm Tr}_{_{(n+1)}}$ and
relations (\ref{bax1}),  one can show

{\bf Lemma~3.1}~~ {\it For all
$X \in \hat{H}_{n}$ and $\forall x,z,$ the following identity is true:
\be
\lb{mapp1}
{\rm Tr}_{_{(n+1)}} \Bigl( T_n(x) \cdot X \cdot T_n(z) \Bigr)  =
(1-x) \, (1-z) \,
{\rm Tr}_{_{(n)}} (X) + \lambda \left(  1 - p \, x \, z \right) \, X \;  ,
\ee
where $T_n(x)$ are Baxterized elements  (\ref{baxtH}) and
$$ p = 1-\lambda D^{(0)} = 1 - (q-q^{-1}) {\rm Tr}_{_{(n+1)}}(1) \; .$$
}
\qed

{}From eq. (\ref{mapp1}), for $p \, x \, z = 1$, we obtain the
"crossing-symmetry relation"
 \be \lb{mapp2}
  \Tr _{_{(n+1)}} \Bigl(
T_n(x) \cdot X \cdot T_n \left( 1/(p x) \right) \Bigr)  =
 \frac{1}{F_p(x)} \, \Tr _{_{(n)}} (X) \;  ,
\ee
where $F_p(x) =\frac{p \, x }{(1-x)(p\, x-1)}$.

\vspace{0.2cm}

\noindent
{\bf Proposition~3.2} (see also \cite{IsOg1}, \cite{Isa07}). {\it Let $y_{n}(x) \in \hat{H}_n$ be any
solution of the RE
(\ref{reflH}). The operators
\be
\lb{tau11}
\tau_{n-1}(x) = \Tr _{_{(n)}} \left( y_{n}(x) \right) \in \hat{H}_{n-1} \; ,
\ee
 form a commutative family of operators
\be
\lb{comfa}
\Bigl[ \tau_{n-1}(x) \, , \; \tau_{n-1}(z) \Bigr]=0 \;\;\; (\forall x,z) \; ,
\ee
in the subalgebra $\hat{H}_{n-1} \subset \hat{H}_{M+1}$.
}

\noindent
{\bf Proof.} Using (\ref{map}), (\ref{mapp2}) and (\ref{reflH}) we find
$$
\begin{array}{c}
\tau_{n-1}(x) \, \tau_{n-1}(z) = Tr_{_{(n)}} \left( y_{n}(x) \, \tau_{n-1}(z) \right) = \\ [0.2cm]
= F_p(x \, z) \,
 \Tr _{_{(n)}}  \left( y_{n}(x) \, \Tr _{_{(n+1)}} \left( T_n (xz) \, y_{n}(z) T_n ((pxz)^{-1}) \right)
\right)=  \\ [0.2cm]
= F_p(x \, z) \,
 \Tr _{_{(n)}} \Tr _{_{(n+1)}} \left( T_n^{-1}(x/z) \, y_{n}(x) \, T_n (xz) \, y_{n}(z)
 \, T_n(x/z) \, T_n ((pxz)^{-1})
\right)= \\ [0.2cm]
= F_p(x \, z) \,
 \Tr _{_{(n)}} \Tr _{_{(n+1)}} \left(  y_{n}(z) \, T_n (xz) \, y_{n}(x)
 \,  T_n ((pxz)^{-1})  \right)= \\ [0.2cm]
=
 \Tr _{_{(n)}}\left(  y_{n}(z) \, \tau_{n-1}(x)  \right) = \tau_{n-1}(z) \, \tau_{n-1}(x) \; ,
\end{array}
$$
where $F_p(x)$ is defined in (\ref{mapp2}). \\
\qed

\vspace{0.3cm}

Now we consider the operators (\ref{tau11}), where solution $y_{n}(x)$ of the reflection equation
is taken in the form (\ref{xxz55}):
 \be
 \lb{bethe02}
 \tau_n(x;\vec{\xi}_{(n)}) = \Tr _{_{(n+1)}} \left( y_{n+1}(x;\vec{\xi}_{(n)})\right) \in \hat{H}_{n}
 \ee
  We stress that the elements (\ref{bethe02})
   are nothing but the analogs of  Sklyanin's transfer-matrices \cite{Skl}
 and the coefficients in the expansion of $\tau_n(x;\vec{\xi}_{(n)})$
 over the variable $x$ (for the homogeneous case $\xi_k=1$) are the Hamiltonians for the
 open Hecke chain models with nontrivial boundary conditions
 which was considered e.g. in \cite{Isa07}.
 Consider this expansion of $\tau_n(x;\vec{\xi}_{(n)})$ for inhomogeneous case:
 \be
 \lb{bethe01}
 \tau_n(x;\vec{\xi}_{(n)})  = \sum_{k=-\infty}^\infty \Phi_k(\vec{\xi}_{(n)})  \, x^k \;  \in \hat{H}_{n} \; .
 \ee
 According to the
 Proposition 3.2, for fixed parameters
 $\vec{\xi}_{(n)}=(\xi_{1},\dots,\xi_n)$, the elements $\Phi_k(\vec{\xi}_{(n)})$
   generate a commutative subalgebra
 $\hat{\cal B}_n(\vec{\xi}_{(n)}) \subset \hat{H}_{n}$. These elements are interpreted as
 Hamiltonians for the inhomogeneous
 open Hecke chain models. Following \cite{MTV}
 we call the subalgebras $\hat{\cal B}_n(\vec{\xi}_{(n)})$ as Bethe subalgebras of the affine
 Hecke algebra $\hat{H}_{n}$.

 First we obtain more explicit form for the generating function of the elements
 $\Phi_k(\vec{\xi}_{(n)}) \in \hat{\cal B}_n(\vec{\xi}_{(n)})$. For this  we substitute
the solution $y_{n+1}(x;\vec{\xi}_{(n)})$ of the reflection equation in the form (\ref{xxz55})
to the transfer-matrix operators (\ref{bethe02}).
Using relation (\ref{mapp1}) we obtain
\be
\lb{tau2}
\begin{array}{c}
\tau_{n}(x;\vec{\xi}_{(n)}) = \Tr _{_{\!\! (n+1)}} \! \left( \!
T_{n}(\frac{x}{\xi_{n}}) \cdots T_2(\frac{x}{\xi_2}) T_1(\frac{x}{\xi_1}) \, y_1(x) \,
T_{1}(x \xi_1) T_2(x \xi_2) \cdots T_{n}(x \xi_{n})\! \right) = \\ [0.2cm]
= (\xi_n -x)(\xi_n^{-1}-x) \, \tau_{n-1}(x;\vec{\xi}_{(n-1)}) + \lambda (1-p \, x^2)
 \, y_{n}(x;\vec{\xi}_{(n-1)}) = \\ [0.2cm]
= \prod\limits_{k=n-1}^n(\xi_k -x)(\xi_k^{-1}-x) \, \tau_{n-2}(x;\vec{\xi}_{(n-2)}) + \\ [0.2cm]
 + \lambda (1-p \, x^2) \left(  (\xi_n -x)(\xi_n^{-1}-x) y_{n-1}(x;\vec{\xi}_{(n-2)})
 + y_{n}(x;\vec{\xi}_{(n-1)}) \right) = \dots = \\ [0.2cm]
= \left( \prod\limits_{k=1}^n(\xi_k -x)(\xi_k^{-1}-x) \right) \tau_{0}(x)
 + \lambda (1-p \, x^2)  \, J_n(x;\vec{\xi}_{(n)}) \; ,
\end{array}
\ee
where the element $\tau_{0}(x) = \Tr _{_{\!\! (1)}} \left(y_1(x) \right) \in \hat{H}_0$
by definition is the central element in $\hat{H}_n$. In equation (\ref{tau2})
we have introduced the notation $J_n(x;\vec{\xi}_{(n)})$ for new explicit generating function
 of the commutative elements $\Phi_k(\vec{\xi}_{(n)}) \in \hat{H}_{n}$:
\be
\lb{tau3}
\begin{array}{c}
J_n(x;\vec{\xi}_{(n)}) = \sum\limits_{r=1}^{n}
   d^n_r(x;\vec{\xi}\;) \;
 y_{r}(x;\vec{\xi}_{(r-1)}) = \\ [0.2cm]
 = \sum\limits_{r=1}^{n}
   d^n_r(x;\vec{\xi}\;) \;
 T_{r-1}(\frac{x}{\xi_{r-1}}) \cdots  T_1(\frac{x}{\xi_1}) \, y_1(x) \,
T_{1}(x \xi_1)  \cdots T_{r-1}(x \xi_{r-1}) \; .
\end{array}
\ee
where we have used the concise notation $d^n_r(x;\vec{\xi}\;)$ for coefficient functions
 \be
 \lb{tau3r}
 d^n_r(x;\vec{\xi}_{(n)}) = \prod\limits_{k=r+1}^n (x-\xi_k)(x-\xi_k^{-1})
 = \prod\limits_{k=r+1}^n \left(1-   \rho_k \, x +x^2 \right) \; ,
 \ee
 \be
 \lb{tau3rr}
 \rho_k = (\xi_k + \xi_k^{-1}) \; .
 \ee

\vspace{0.2cm}

\noindent
{\bf Remark.}  From Yang-Baxter equation (\ref{ybeH}) and
 reflection equation (\ref{reflH}) we deduce
 \be
 \lb{symT}
 \begin{array}{c}
\tau_{n}(x;\vec{\xi}_{(n)}) \, T_k(\xi_{k+1}/\xi_k)   = T_k(\xi_{k+1}/\xi_k) \,
 \tau_{n}(x; {\bf s}_k \cdot \vec{\xi}_{(n)}) \; , \;\;\; k=1,\dots,n-1 \; , \\ [0.2cm]
\tau_{n}(x;\vec{\xi}_{(n)}) \, \bar{y}_k(\xi_{k};\vec{\xi}_{(k-1)})   = \bar{y}_k(\xi_{k};\vec{\xi}_{(k-1)})  \,
 \tau_{n}(x; {\bf I}_k \cdot \vec{\xi}_{(n)}) \; , \;\;\; k=1,\dots,n \; ,
 \end{array}
 \ee
 where ${\bf s}_k \cdot \vec{\xi}_{(n)} \equiv  (\xi_1, \dots, \xi_{k-1},\xi_{k+1},\xi_k,\xi_{k+2}, \dots \xi_n)$, i.e.
 ${\bf s}_k$ is the transposition of two parameters $\xi_k$ and $\xi_{k+1}$,
 and $ {\bf I}_k \cdot \vec{\xi}_{(n)} \equiv
  (\xi_1, \dots, \xi_{k-1},\xi^{-1}_k,\xi_{k+1}, \dots \xi_n)$. It means that
 the Bethe subalgebras generated by
 transfer-matrix type elements $\tau_{n}(x;\vec{\xi}_{(n)})$, $\tau_{n}(x; {\bf s}_k \cdot \vec{\xi}_{(n)})$
  and $\tau_{n}(x; {\bf I}_k \cdot \vec{\xi}_{(n)})$ are equivalent. It is clear that the
  symmetry (\ref{symT}) is also valid for the generating functions (\ref{tau3}).

 \section{Bethe subalgebra for the Hecke algebra}

The Hecke algebra $H_n$ is the quotient of the affine Hecke algebra $\hat{H}_n$
by the relation $\hat{y}_1=1$.
Thus, one can obtain the generating function for the elements of Bethe
subalgebra of usual Hecke algebra $H_n$ if we substitute into
(\ref{tau3}) the trivial solution
$y_1(x)=1$ of the reflection equation. Then to simplify the function (\ref{tau3})
for the case of $H_n$ we first transform
the elements $y_{r}(x;\vec{\xi}_{(r-1)})$ given in (\ref{xxz55}). For this
we use relations (\ref{baxtH}) and identities
\be
\lb{tau4}
\begin{array}{c}
 T_k(\frac{x}{\xi_k}) \, T_{k}(x \xi_k) =
\lambda T_k (x^2)  + (x-\xi_k)(x-\xi_k^{-1}) =  \\ [0.2cm]
= \lambda (1-x^2) T_k + [ \lambda^2 x^2 + (x-\xi_k)(x-\xi_k^{-1})]  \; ,
\end{array}
\ee
which can be deduced from (\ref{bax1}). For $y_1(x)=1$, applying (\ref{tau4})
many times, we obtain new representation for the elements (\ref{xxz55}):
\be
\lb{tau3b}
\begin{array}{c}
y_{n+1}(x;\vec{\xi}_{(n)}) =
T_{n}(\frac{x}{\xi_{n}}) \cdots T_2(\frac{x}{\xi_2}) T_1(\frac{x}{\xi_1}) \,
T_{1}(x \xi_1) T_2(x \xi_2) \cdots T_{n}(x \xi_{n}) =  \\ [0.2cm]
 = T_{n}(\frac{x}{\xi_{n}}) \cdots T_2(\frac{x}{\xi_2})
 \left( \lambda (1-x^2) T_1 + [ \lambda^2 x^2 + (x-\xi_1)(x-\xi_1^{-1})]
 \right) T_2(x \xi_2) \cdots T_{n}(x \xi_{n}) =  \\ [0.2cm]
  = \lambda (1-x^2) \; \widetilde{y}_{n+1}(x;\vec{\xi}_{(n)})
+ c_{n+1} (x;\vec{\rho}_{(n)})  \; ,
\end{array}
\ee
where
\be
\lb{tau3y}
\begin{array}{c}
\widetilde{y}_n(x;\vec{\xi}_{(n-1)})  = \sum\limits_{k=1}^{n-1} c_k (x;\vec{\rho}_{(k-1)})
 T_{n-1}(\frac{x}{\xi_{n-1}}) \cdots T_{k+1}(\frac{x}{\xi_{k+1}})
T_{k} T_{k+1}(x \xi_{k+1}) \cdots  T_{n-1}(x \xi_{n-1}) \; ,
\end{array}
 \ee
 parameters $\rho_k$ were defined in (\ref{tau3rr})
and for the coefficient functions $c_k$ we have $c_1=1$,
$$
c_k (x;\vec{\rho}_{(k-1)}) \equiv
\prod\limits_{j=1}^{k-1} \left( 1-x \rho_j  + (1+\lambda^2) x^2 \right)
\;\;\;\;\; (\forall k \geq 2) \; .
$$
Using representation (\ref{tau3b})
it is convenient to redefine the generating function (\ref{tau3}) once again
\be
\lb{tau3J}
\begin{array}{c}
J_n(x;\vec{\xi}_{(n)})
 = \sum\limits_{r=1}^{n}
   d^n_r(x;\vec{\xi}\;) \;
 \left( \lambda (1-x^2) \widetilde{y}_{r}(x;\vec{\xi}_{(r-1)}) + c_{r}(x;\vec{\rho}_{(r-1)})  \right) = \\ [0.2cm]
 = \lambda (1-x^2) \widetilde{J}_n(x;\vec{\xi}_{(n)}) + \sum\limits_{r=1}^{n}
   d^n_r(x;\vec{\xi}\;) \; c_{r}(x;\vec{\rho}_{(r-1)})  \; .
\end{array}
\ee
For new function $\widetilde{J}_n(x;\vec{\xi}_{(n)})$ which generate elements
of the Bethe subalgebra ${\cal B}_n(\vec{\xi}_{(n)}) \subset H_n$ we obtain
the recurrent relations
\be
 \lb{tau4c}
\widetilde{J}_2= \widetilde{y}_2 =  T_1 \; , \;\;\;
 \tilde{J}_n  = (1-x \rho_n + x^2) \tilde{J}_{n-1}  + \widetilde{y}_n \;
= \; \sum\limits_{k=2}^n \, d^n_k(x;\vec{\rho}) \, \widetilde{y}_k \; ,
 \ee
where coefficients $d^n_k(x;\vec{\rho})$ were defined in (\ref{tau3r}). Using
the recurrence relation $(30)$ we can compute the Hamiltonian of our
(integrable) system, namely,
$$ {\partial \over \partial x} \widetilde{J}_n(x;\vec{\xi}_{(n)}) |_{x=0} =
\sum_{1 \le i < j \le n} (\rho_i + \rho_j) T_{(ij)} + \lambda~\sum_{1 \le i < j < k < n} \bigl(\xi_j^{-1} T_{(ij)} T_{(jk)} +\xi_j T_{(jk)} T_{(ij)} \bigr).$$

 At the end of this Section we present the explicit expressions for first few
elements $\widetilde{y}_n$  and $\widetilde{J}_n$ for $n \geq 2$:
 \be
 \lb{tau4d}
\begin{array}{c}
\widetilde{y}_2 = T_1 \; , \;\; \widetilde{y}_3 = T_{2}(\frac{x}{\xi_{2}})
T_{1} T_{2}(x \xi_{2}) + \left( 1-x \rho_1  + (1+\lambda^2) x^2 \right)T_2
\; , \\ [0.3cm]
\widetilde{y}_4 = T_{3}(\frac{x}{\xi_{3}}) T_{2}(\frac{x}{\xi_{2}})
T_{1} T_{2}(x \xi_{2})T_{3}(x \xi_{3}) +
 \left( 1-x \rho_1  + (1+\lambda^2) x^2 \right)T_{3}(\frac{x}{\xi_{3}})
T_{2} T_{3}(x \xi_{3}) \\ [0.2cm]
+  [ 1-x \rho_1  + (1+\lambda^2) x^2][ 1-x \rho_2  + (1+\lambda^2) x^2] T_3 \; , \\ [0.3cm]
\widetilde{y}_5 = T_{4}(\frac{x}{\xi_{4}}) ... T_{2}(\frac{x}{\xi_{2}})
T_{1} T_{2}(x \xi_{2})...T_{4}(x \xi_{4}) +
 c_2(x;\rho_1) \, T_{4}(\frac{x}{\xi_{4}})
 T_{3}(\frac{x}{\xi_{3}}) T_{2} T_{3}(x \xi_{3})T_{4}(x \xi_{4}) + \\ [0.2cm]
  + c_3(x;\rho_1,\rho_2) \,  T_{4}(\frac{x}{\xi_{4}}) T_{3} T_{4}(x \xi_{4})
  + c_4(x;\rho_1,\rho_2,\rho_3) \, T_4 \; .
  \end{array}
 \ee
 \be
 \lb{tau5d}
\begin{array}{c}
\widetilde{J}_2=  T_1 \; , \;\;\; \widetilde{J}_3=  \Bigl( 1-x \rho_3  +
(1+\lambda^2) x^2 \Bigr) T_1 + \Bigl( 1-x \rho_1  +
(1+\lambda^2) x^2 \Bigr) T_2  + \\ [0.2cm]
+ \left( 1-x \rho_2  +  x^2 \right) T_2 T_1 T_2  +
\lambda x   (\xi_2 - x) \, T_2 T_1 + \lambda x  \, (\xi_2^{-1} - x) \,
T_1 T_2 = \\ [0.2cm]
T_1+T_2+T_1 T_2 T_1 - ( \rho_3 T_1+ \rho_2 T_1 T_2 T_1 +\rho_1 T_2 -
\lambda \xi_2~ T_2 T_1 - \lambda \xi_2^{-1} T_1 T_2) x + \\ [0.2cm]
\Bigl((1+\lambda^2)(T_1+T_2+T_1 T_2 T_1)- \lambda (T_1 T_2 +T_2 T_1 +
\lambda~ T_1 T_2 T_1) \Bigr) x^2

\; .
 \end{array}
 \ee
Therefore the Bethe subalgebra ${\cal B}_3(\vec{\xi}_{(3)})$ is generated by
the central elements $C_1= T_1+T_2 + T_{(13)}$~and~$C_2 = T_1 T_2 +T_2 T_1 +
 \lambda T_{(13)},$~and~the element
$$D= \rho_3 T_1 + \rho_2 T_{(13)} +\rho_1 T_2 - \lambda \xi_2^{-1} T_2 T_1 -
\lambda \xi_2~ T_1 T_2.$$
Here we used notation $T_{(13)}:= T_1 T_2 T_1.$ Now let us introduce the
following elements
$$\theta_1 := \theta_1^{\vec{\xi}_{(3)}}={D-\rho_1 C_1 \over (\rho_2-\rho_1)(\rho_3 -\rho_1)} ={T_1 \over
\rho_2-\rho_1} +
{T_{(13)} \over \rho_3-\rho_1} -{\lambda \xi_2 ~T_2 T_1 + \lambda \xi_2^{-1} T_1 T_2  \over (\rho_2-\rho_1)(\rho_3-\rho_1)},$$
$$\theta_2 :=  \theta_1^{\vec{\xi}_{(3)}}={D -\rho_2 C_1 \over (\rho_1-\rho_2)(\rho_3-\rho_2)} =
{T_1 \over \rho_1-\rho_2}+
{T_2 \over \rho_3-\rho_2} -{\lambda \xi_2~ T_2 T_1 + \lambda \xi_2^{-1} T_1 T_2
\over (\rho_1-\rho_2)(\rho_3-\rho_2)},$$
$$\theta_3 := \theta_1^{\vec{\xi}_{(3)}}={D- \rho_3 C_1 \over (\rho_1-\rho_3)(\rho_2-\rho_3)} =
{T_2 \over \rho_2-\rho_3} +
{T_{(13)} \over \rho_1-\rho_3} -{\lambda \xi_2~ T_2 T_1 + \lambda \xi_2^{-1}
T_1 T_2 \over (\rho_2-\rho_3)(\rho_1-\rho_3)}.$$
One can check that
$$ \theta_1+\theta_2+\theta_3=0,~~\rho_1 \theta_1 +\rho_2 \theta_2 +\rho_3 \theta_3 = C_1,~~ (\lambda^2+3) C_2 =C_1^2 -2 \lambda C_1 -3,$$
and the elements ~~$\theta_1^{\vec{\xi}_{(3)}},\theta_2^{\vec{\xi}_{(3)}},
\theta_3^{\vec{\xi}_{(3)}}$~~pairwise commute and generate the Bethe
subalgebra ~ ${\cal B}_3(\vec{\xi}_{(3)})$. ~~Our goal is to show that a
similar set of generators exist for the Bethe algebra~~
${\cal B}_n(\vec{\xi}_{(n)})$~~ for arbitrary $n$.

 To state our main result
of this Section we need to introduce a bit of notation. First of all, for a
pair of integers $1 \le i < j \le n$~~let us introduce the elements ~$T_{(ij)}:=
T_{j-1} \cdots T_{i+1} T_{i} T_{i+1} \cdots T_{j-1},~~1 \le i < j \le n.$ Now
 let $B \subset [1,2,\ldots,n]$ be a subset, define inductively the elements
~$T(B):= T^{{\vec{\xi}}_{(n)}}(B)$ as follows

$\bullet$  ~~$ T(\{ b \})=0,$~~~$T(\{ a < b \}) =T_{(ab)},$

$\bullet$~~$T(\{a < b < c  < \ldots < d \}) = \xi_a~ T(\{b < c  < \ldots < d \})~ T_{(ab)}
+ \xi_a^{-1} T_{(ab)}~ T(\{ b < c  < \ldots < d \}).$ \\
For example, \\
$T(\{a < b < c < d \}) = \xi_a\xi_b~ T_{(cd)} T_{(bc)} T_{(ab)}+
\xi_a \xi_b^{-1} T_{(bc)} T_{(cd)} T_{(ab)}+ \xi_a^{-1} \xi_b~ T_{(ab)}
 T_{(cd)} T_{(bc)} + \xi_a^{-1} \xi_b^{-1} T_{(ab)} T_{(bc)} T_{(cd)}.$

Using the notation introduced above, let us define the following elements
$$\theta_a^{\vec{\xi}_{(n)}} = \sum_{ B \subset [1,\ldots, n] \atop a \in B}
\lambda^{ |B|-2} ~~{ T(B) \over \prod_{ b \in B \atop b \not= a} (\rho_a - \rho_b)},~~~a= 1,\ldots,n.$$

\begin{thm} ${}$

 The elements ~~$ \theta_a^{\vec{\xi}_{(n)}},$~$ a=1, \ldots, n$ ~~ mutually
commute,  generate the Bethe subalgebra ~${\cal B}_n(\vec{\xi}_{(n)})$~ and satisfy the following properties

$\bullet$~~$\theta_1^{\vec{\xi}_{(n)}}+ \ldots + \theta_n^{\vec{\xi}_{(n)}} =0,$~~$\sum_{j=1}^{n} \xi_j \theta_j^{\vec{\xi}_{(n)}} = \sum_{1 \le i < j \le n}
T_{(ij)},$

$\bullet$~~the elementary symmetric polynomials $e_j(\theta_1^{\vec{\xi}_{(n)}}, \ldots, \theta_n^{\vec{\xi}_{(n)}}),$~$j=2,\ldots,n,$~ generate the center
of the Hecke algebra $H_n,$
\end{thm}
Clearly,
$$\theta_a^{\vec{\xi}_{(n)}} = \sum_{b \not= a} {T_{(ab)} \over \rho_a - \rho_b} ~+ \lambda \bigl( \ldots \bigr) .$$
In other words the elements $\{\theta_a^{\vec{\xi}_{(n)}},~a=1,\ldots,n \}$
are a {\it lift} of the Gaudin elements $\{ g_a ({\vec{\xi}_{(n)}}) := \sum_{j \not= a}~~(\rho_a-\rho_j)^{-1}~s_{aj},~a=1,\ldots,n \}$ from the group algebra
$\C [\S_{n}]$ of the symmetric group $\S_{n}$ to the Hecke algebra $H_n \otimes \C .$

Using
the recurrence relation $(30)$ we can compute the Hamiltonian
 ${\cal{H}}^{\vec{\xi}_{(n)}}$ of our
(integrable) model, namely,
$$ { \cal H}^{\vec{\xi}_{(n)}} = {\partial \over \partial x} \widetilde{J}_n(x;\vec{\xi}_{(n)}) |_{x=0} =
\sum_{1 \le i < j \le n} (\rho_i + \rho_j) T_{(ij)} +
 \sum_{B \subset [1,\ldots,n]  \atop |B| \ge 3} ~(- \lambda)^{|B|-2}~T(B).$$

~Finally we remark that
the  example above shows that the set of all
Bethe's subalgebras in the Hecke algebra $H_n$ does not coincide with the set
of all maximal commutative subalgebras in $H_n,$ if $ n \ge 3.$

 \section{Symmetric group limit and Gaudin model.}

Let us consider the special classical limit when $q \to 1$ while parameters $x$ and
$\xi_k$ are fixed.
For $q=1$ or $\lambda=0$ in view of (\ref{ahecke}), (\ref{tau4}) the
Hecke algebra $H_{M+1}$ is degenerated to the symmetric group algebra ${\cal S}_{M+1}$, i.e.
$T_k = T^{-1}_k = s_{k,k+1} =s_k$ are elementary transpositions of $k$ and $k+1$. In this limit we
 have $T_{k}(x) = (1-x)s_{k}$ and formulas (\ref{tau3y}) and (\ref{tau4d}) are simplified
\be
 \lb{tau4e}
\begin{array}{c}
\widetilde{y}_2 = s_1 \; , \;\; \widetilde{y}_3 = [ 1-x \rho_2  + x^2 ] s_{2}
 s_{1} s_{2} + [ 1-x \rho_1  + x^2 ] s_2 \; , \\ [0.2cm]
\widetilde{y}_4 = [ 1-x \rho_2  +  x^2][ 1-x \rho_3  +  x^2] s_{3} s_{2}
s_{1} s_{2}s_{3} +
 [ 1-x \rho_1  +  x^2 ] [ 1-x \rho_3  +  x^2]s_{3}
s_{2} s_{3} \\ [0.2cm]
+  [ 1-x \rho_1  +  x^2][ 1-x \rho_2  +  x^2] s_3 \; , \; , \dots ,
 \end{array}
 \ee
 \be
 \lb{tau4f}
\begin{array}{c}
  \widetilde{y}_k = \left( \prod\limits_{m=1}^{k-1} [ 1-x \rho_k  +  x^2] \right)
   \sum\limits_{j=1}^{k-1} \frac{s_{j,k}}{[ 1-x \rho_j  +  x^2]}   \; ,
\end{array}
 \ee
 where $s_{j,k} = s_{k-1} ...s_{j+1} s_{j} s_{j+1}...s_{k-1}$ are
 transpositions in ${\cal S}_{M+1}$, i.e. $s_{j,k} = s_{k,j}$.
 Substitution of (\ref{tau4f}) into (\ref{tau4c}) gives
 $$
 \begin{array}{c}
  \widetilde{J}_n(x;\vec{\rho}_{(n)}) = \; \sum\limits_{k=2}^n \,
   \prod\limits_{m=k+1}^n \left( 1-x \rho_m  +  x^2 \right)  \, \widetilde{y}_k = \\ [0.2cm]
 =   \sum\limits_{k=2}^n \, \sum\limits_{j=1}^{k-1} \,
  \prod\limits_{\stackrel{m=1}{m\neq j,k}}^n \left( 1-x \rho_m  +  x^2 \right)
     \, s_{j,k}
   =  \prod\limits_{m=1}^n \left( 1-x \rho_m  +  x^2 \right)
 \sum\limits_{k>j}^n \,
   \frac{s_{j,k} }{  [ 1-x \rho_j  +  x^2][ 1-x \rho_k  +  x^2]}  \; ,
  \end{array}
  $$
After the renormalization
 $$
 \widetilde{J}_n(x;\vec{\rho}_{(n)}) \to  \widetilde{J}'_n(x;\vec{\rho}_{(n)}) =
 \frac{x^2 \; \widetilde{J}_n(x;\vec{\rho}_{(n)})}{\prod\limits_{m=1}^n \left( 1-x \rho_m  +  x^2 \right)}
 $$
and change of variables $u =  x+1/x$ we obtain the generating function
for Bethe subalgebra of symmetric group algebra ${\cal S}_n$ in the form
    \be
 \lb{tau4h}
  \begin{array}{c}
 \widetilde{J}'_n(x;\vec{\rho}_{(n)}) =
 \frac{x^2}{2} \, \sum\limits_{\stackrel{k,j=1}{k \neq j}}^n \,
   \frac{ s_{j,k} }{  [ 1-x \rho_j  +  x^2][ 1-x \rho_k  +  x^2]} =
   \frac{1}{2} \, \sum\limits_{\stackrel{k,j=1}{k \neq j}}^n \,
   \frac{ s_{j,k} }{[ u - \rho_j][u- \rho_k ]} =
   \frac{1}{2} \, \sum\limits_{\stackrel{k,j=1}{k \neq j}}^n \,
   \frac{ s_{j,k} }{\rho_j -\rho_k}
   \frac{\rho_j -\rho_k}{[ u - \rho_j][u- \rho_k ]} = \\ [0.2cm]
  =  \frac{1}{2} \,  \sum\limits_{\stackrel{k,j=1}{k \neq j}}^n \,
   \frac{ s_{j,k} }{\rho_j -\rho_k}
   \left( \frac{1}{[ u - \rho_j]} - \frac{1}{[u- \rho_k ]} \right) =
   \sum\limits_{\stackrel{k,j=1}{k \neq j}}^n \,
   \frac{ s_{j,k} }{\rho_j -\rho_k}
  \frac{1}{[ u - \rho_j]}  \, .
   \end{array}
\ee
The commuting Hamiltonians $H^{[n]}_j$ for Gaudin model are obtained from
(\ref{tau4h}) as residues for $u \to \rho_j$:
 $$
 H^{[n]}_j = \left. {\rm res}(\widetilde{J}'_n(x;\vec{\rho}_{(n)})) \right|_{u=\rho_j} =
 \sum\limits_{\stackrel{k=1}{k \neq j}}^n \,
   \frac{ s_{j,k} }{\rho_j -\rho_k} \; .
 $$

 The right hand side of (\ref{tau4h}), after the change of variables $\rho_k \to z_k$,
  can be represented in the form
 $$
  \frac{1}{\prod\limits_m [ u - z_m]}
  \left( \sum\limits_{\stackrel{k,j=1}{k \neq j}}^n \, \frac{ s_{j,k} }{z_j -z_k}
  \prod\limits_{\stackrel{m}{m \neq j}} [ u - z_m] \right) \; ,
 $$
and the expression in the brackets is related
(up to the shift by a scalar function) to the generating
function $\Phi^{[n]}_2(u)$ of the Bethe subalgebra which was presented in
\cite{MTV} (see Remark after the Theorem 4.3 in \cite{MTV}).

\vspace{0.5cm}

\noindent

{\bf Remark 1.} There is another semi-classical limit $q \to 1$
for the constructions considered above which is called \underline{Yangian limit}.
In this case we have to consider substitution
 $$
 x = q^{-2 u} \; , \;\;\; \xi_j = q^{-2 z_j}  \; , \;\;\; q = e^h \; ,
 $$
 and take the limit $ h \to 0$. Then we have
 \be
 \lb{spin}
 \begin{array}{c}
x = 1 - 2 h u + ... \; , \;\;\; \lambda = 2h + ...
 \; , \;\;\; T_k(x) = 2 h ( u s_k + 1) + ...  \; , \\ [0.2cm]
 T_k = s_k + h \, s_k' + \dots \; , \;\;\;   T_k^2 = 1 + 2 \, h \, s_k + \dots \; , \;\;\;
   s_k s'_k + s'_k  s_k = 2 s_k \; , \\ [0.2cm]
 (1 - \rho_j x + x^2) = (1 - \xi_j x) (1 - {x \over \xi_j}) = 4 h^2 (u^2 - z_j^2) + ... \; , \\ [0.2cm]
 (1 - \rho_j x + (1+\lambda^2) x^2) = 4 h^2 (u^2 - z_j^2 +1) + ... \; , \\ [0.2cm]
 \hat{y}_1 = 1 + 2h \, \ell_1  + \dots \; , \;\;\;  \hat{y}_k = 1 + 2 h \, \ell_k' + \dots \; , \\ [0.2cm]
  T_{k}(\frac{x}{\xi_{k}}) =  2 h ( (u - z_k) s_k + 1) \; , \;\;\;
  T_{k}(x \xi_{k}) =  2 h ( (u+z_k) s_k + 1) \; ,
 \end{array}
 \ee
where $s_k \in \mathbb{S}_{M+1}$ are elementary transpositions, dots $\dots$ means the
higher orders in $h$ and we have used notation
\be
 \lb{bax5}
\ell'_1 = \ell_1 \; , \;\;\;  \ell'_{k+1} = s_{k} + s_k \ell_{k}' s_k \; ,
 \ee
for commutative set of elements $\{ \ell_k' \}$. Moreover it is easy to see that
~~$s_k \ell_{k}'= \ell_{k+1}' s_{k} -1, $~$k=1,\ldots, M,$~and the elements
$\{s_k , \ell_{1}',\ell_{k+1}',~k=1,\ldots,M \}$ generate the degenerate affine Hecke
algebra \cite{Drinf}.

We call the limit (\ref{spin}) as Yangian limit since
for the modified Baxterized element ${\widetilde{T}}_k(x)$ given in
(\ref{bax3}) we obtain (see the first line in (\ref{spin})):
 \be
 \lb{bax4}
{\widetilde{T}}_k(x) = \frac{u s_k + 1}{u+1} + h \, \dots \; .
 \ee
In the matrix representation $\rho$
of the symmetric group $\mathbb{S}_{M+1}$ we have $\rho(s_k) = P_{kk+1}$, where $P_{kk+1}$ are permutation matrices
which act in the products of $n$-dimensional vector spaces
$V^{\otimes (M+1)}$ and permute factors with numbers $k$ and $k+1$. In this
representation the Baxterized element (\ref{bax4}) is nothing but Yangian
$R$-matrix which is used for the definition of the $g\ell(n)$-type Yangian.

\vspace{0.5cm}

\noindent
{\bf Remark 2.} Now we show that the elements $\ell_1$ commute as following:
\be
 \lb{gl04}
[\ell_1,\ell_2] = s_1 \ell_1 - \ell_1 s_1 \; ,
 \ee
where we denoted $\ell_2 := s_1 \, \ell_1 \, s_1$. It follows from the relation $[\ell_1', \, \ell'_2]=0$, but
we check this identity directly. First, we start with the relation
 (see (\ref{afheck})) in the affine Hecke algebra
 \be
 \lb{gl01}
 T_1 \, \hat{y}_1 \, T_1 \,  \hat{y}_1 = \hat{y}_1 \, T_1 \, \hat{y}_1 \, T_1   \; ,
 \ee
 and change here the affine generators $\hat{y}_1 \to K_1$:
 \be
 \lb{gl02}
 \hat{y}_1 = 1 + (q-q^{-1}) K_1  \; .
 \ee
 The substitution of (\ref{gl02}) into (\ref{gl01})
 and using of (\ref{ahecke}) gives for $(q-q^{-1}) \neq 0$:
 \be
 \lb{gl03}
 T_1 K_1 T_1 K_1 -   K_1 T_1 K_1 T_1 = K_1 T_1 - T_1 K_1 \; .
 \ee
 Now we take here the Yangian limit (\ref{spin}), i.e. we put
 $K_1 = \ell_1 + h \dots$ and $T_1 = s_1 + h \dots$. Then, in the zero order
 approximation in $h$, eq. (\ref{gl03}) is converted to (\ref{gl04}).
 In the matrix representation (see the end of the Remark 1.) when $\rho(s_k) = P_{kk+1}$,
 $$
 \rho(\ell_1) = L \otimes \underbrace{1 \otimes \cdots \otimes 1}_{M} \in End(V^{\otimes (M+1)})
 $$
  and $L=||L_{ij}||_{i,j=1,...,n}$ is a $n \times n$ matrix with noncommutative entries,
   relation (\ref{gl04}) gives the defining relations for
  generators $L_{ij}$ of $g\ell(n)$ Lie algebra.

\vspace{0.5cm}

\noindent
{\bf Remark 3.} The using of the Yangian limit (\ref{spin}) in (\ref{tau3})
will give the generating function for generators of Bethe subalgebras for
degenerate affine Hecke algebra. Here we present first few cases for $n=1,2,3$:
$$
J_1(x,\vec{\xi}_{(1)}) = y_1(x)  = \frac{\y_1 - \xi \, x }{\y_{1} - \xi \, x^{-1}}
 \;\; \stackrel{h \to 0}{\longrightarrow} \;\; \frac{\ell_1 + (z+u) }{\ell_{1} + (z-u) } + h \dots \; ,
$$
$$
J_2(x,\vec{\xi}_{(2)}) = (x^2 - \rho_2 x +1) \, y_1(x) + T_1( x/\xi_1) \, y_1(x) \, T_1( x \xi_1)
  \;\; \stackrel{h \to 0}{\longrightarrow}
$$
$$
 4 h^2 \left( (u^2 - z_2^2)  \frac{\ell_1 + (z+u) }{\ell_{1} + (z-u) }  +
((u-z_1) s_1 +1) \frac{\ell_1 + (z+u) }{\ell_{1} + (z-u) } ((u+z_1) s_1 +1) \right) + h^3 \dots \; ,
$$
$$
 \begin{array}{c}
J_3(x,\vec{\xi}_{(3)}) = (x^2 - \rho_3 x +1) \Bigl( (x^2 - \rho_2 x +1) \, y_1(x) +
  T_1( x/\xi_1) \, y_1(x) \, T_1( x \xi_1) \Bigr) +   \\ [0.2cm]
  + T_2( x/\xi_2) \, T_1( x/\xi_1) \, y_1(x) \, T_1( x \xi_1) \, T_2( x \xi_2)
  \;\; \stackrel{h \to 0}{\longrightarrow}
  \end{array}
$$
$$
 \begin{array}{c}
 \displaystyle{ (2 h)^4 \left( (u^2 - z_3^2)  (u^2 - z_2^2)  \frac{\ell_1 + (z+u) }{\ell_{1} + (z-u) }  + (u^2 - z_3^2)
((u-z_1) s_1 +1) \frac{\ell_1 + (z+u) }{\ell_{1} + (z-u) } ((u+z_1) s_1 +1) \right. + } \\ [0.4cm]
 \displaystyle{  \left.
 + ((u-z_2) s_2 +1)  ((u-z_1) s_1 +1) \frac{\ell_1 + (z+u) }{\ell_{1} + (z-u) }
  ((u+z_1) s_1 +1) ((u+ z_2) s_2 +1) \right)   +  h^5 \dots }  \; ,
  \end{array}
$$
where we have used (\ref{soluH}) and taken $\xi = q^{-2z}$. The expansion of elements
$\left. J_n(x,\vec{\xi}_{(n)}) \right|_{h \to 0}$ over the spectral parameter
$u$ gives the set of commuting elements $\Phi_k(z_1,\dots,z_n)$ (cf. eq. (\ref{bethe01})):
 $$
 \left. J_n(x,\vec{\xi}_{(n)}) \right|_{h \to 0} =
  (2 h)^{2(n-1)} \sum_{k=0}^\infty \Phi_k(z_1,\dots,z_n) \, u^k + h^{2n-1} \dots \; .
 $$
 The elements $\Phi_k(z_1,\dots,z_n)$ are generators of the
 Bethe subalgebras in the degenerate affine Hecke algebra $\{s_m , \ell_{1},~m=1,\ldots,n-1 \}$
 and $\Phi_k(z_1,\dots,z_n)$ can be interpreted
 as Hamiltonians for inhomogeneous open spin chains. In particular for homogeneous case
 $z_m =0$, $z \neq 0$ of such open spin chain we obtain the local Hamiltonian ${\cal H}(z)$:
 $$
 \left. \partial_u \; J_n \Bigl(q^{-2u},(1,...,1) \Bigr) \right|_{h \to 0, u=0}  \; \sim \;
 h^{2(n-1)} \, {\cal H}(z) + h^{2n-1} \dots \;\; , \;\;\;\;\;
 {\cal H}(z)  = \left( \sum_{m=1}^{n-1} s_m + \frac{1}{\ell_1 + z} \right) \; .
 $$

\end{document}